\documentclass[a4paper,oneside,11pt,reqno]{amsart}
\pdfoutput=1
\usepackage{amsmath, latexsym, amsfonts, amssymb, amsthm, amscd, stmaryrd,setspace,bm,graphicx}
\usepackage[utf8]{inputenc}

 \theoremstyle{plain}
 \newtheorem{theorem}{Theorem}
 \newtheorem*{theorem*}{Theorem}
 \newtheorem{proposition}{Proposition}[section]
 \newtheorem*{proposition*}{Proposition}
 \newtheorem{lemma}{Lemma}[section]

 \newtheorem{rem}{Remark}[section]
  \newtheorem*{rem*}{Remark}
 
 \newtheorem*{definition*}{Definition}
 
\newtheorem*{fact*}{Fact}

 \renewcommand{\tilde}{\widetilde}          
 \DeclareMathSymbol{\leqslant}{\mathalpha}{AMSa}{"36} 
 \DeclareMathSymbol{\geqslant}{\mathalpha}{AMSa}{"3E} 
 \DeclareMathSymbol{\eset}{\mathalpha}{AMSb}{"3F}     
 \renewcommand{\leq}{\;\leqslant\;}                   
 \renewcommand{\geq}{\;\geqslant\;}                   

 \newcommand{\R}{\mathbb{R}}
 \newcommand{\Z}{\mathbb{Z}}
 \newcommand{\N}{\mathbb{N}}
 
\newcommand{\eql}{\stackrel{(law)}{=}}
 \newcommand{\ind}{\textbf{1}}
\newcommand{\yb}{\bar{y}}

\newcommand{\IP}{\mathbb{P}}
\newcommand{\E}{\mathbb{E}}

\newcommand{\8}{\infty}

\newcommand{\tildeP}{\tilde{P}}
\newcommand{\bP}{\textbf{P}}
\newtheorem{remark}[theorem]{Remark}
\title[Arcsine law and internal DLA generated by Sinai's walk]{The Arcsine law as the limit of the internal DLA cluster generated by Sinai's walk}
\author{N. Enriquez$^{1,*}$}
\thanks{$*$ Modélisation aléatoire de Paris 10 (MODAL'X)-Université Paris Ouest Nanterre La Défense}
\thanks{\textit{email $1$:} nenriquez@u-paris10.fr}
\thanks{\textit{email $2$:} clucas@clipper.ens.fr}
\thanks{\textit{email $3$:} fsimenhaus@u-paris10.fr}
\author{C. Lucas$^{2,*}$}
\author{F. Simenhaus$^{3,*}$}

\begin{document}
\begin{abstract}
 We identify the limit of the internal DLA cluster generated by Sinai's walk as the law of a functional of a Brownian motion which turns out to be a new interpretation of the Arcsine law.
\end{abstract}
\maketitle
 \normalsize

\section{Introduction}
The internal diffusion limited aggregation method was first introduced by Diaconis and Fulton in 1982 (see \cite{DIA_FUL}) and gives a protocol for building a sequence of growing random sets $A(n),\ n\geq 0$ using random walks. At each iteration, the set $A(n+1)$ is obtained from $A(n)$ by addition of the first site visited by a walk starting from $0$ outside $A(n)$.
When the walk is a simple random walk on $\Z^d$, the cluster has the ball as asymptotic limit shape (\cite{LBG}). In the special and trivial case of dimension one, the cluster $A(n)$ is an interval denoted by $[g_n,d_n]$ and $d_n/n$ converges almost surely to $1/2$. In this paper, we consider the case where the cluster is generated by one dimensional walks evolving in an inhomogeneous random medium. More precisely we are dealing with one dimensional recurrent random walks in random environment often called Sinai's walk. One can rapidly realize that depending on the profile of the potential associated to the environment the cluster may be drastically asymmetric. One can even prove that $\limsup d_n/n=1$ and $\liminf d_n/n=0$ almost surely (Theorem $\ref{th:quenched}$). But beyond this rough result we are able to prove, under the annealed measure, a convergence in law for $d_n/n$ towards the Arcsine law (Theorem $\ref{th:mainth}$). The Arcsine law comes out from a new functional of the Brownian motion which in our case is the scaling limit of the potential. This functional stems from the exploration of the potential by the growing cluster. It involves in Sinai's terminology the largest valley of \textit{width smaller than one} containing the origin. Let us remind that the limit law of the walk in Sinai's theorem after scaling by $\ln^2 n$ involves the smallest valley of \textit{height one}.

Finally, we mention that this result provides a new understanding of the Arcsine law which is different from the two classical ones (last zero of the Brownian motion before $1$, and the time spent by the Brownian motion above $0$ before $1$).

\section{Notations and main results}

We first introduce Sinai's random walk in random environment. An environment $\omega$ is a collection $(\omega(i))_{i \in \Z}$ of numbers in $[0,1]$. We denote by $\Omega := [0,1]^{\Z}$ the set of environments.

For a given environment $\omega$, we define a Markov chain $(X_n)_{n \in \N}$ of law $\tildeP_{\omega}$, often called \textit{quenched law}, by
\begin{eqnarray*}
\tildeP_{\omega}(X_0 = 0) &=& 1, \\
\end{eqnarray*}
and for every $x\in\Z$ and $n\in\N$,
\begin{eqnarray*}
\tildeP_{\omega}(X_{n+1}  = x+1 | X_n =x) & = & \omega(x),  \\
\tildeP_{\omega}(X_{n+1}  = x-1 | X_n =x) & = & 1 - \omega(x).
\end{eqnarray*}
We endow $\Omega$ with its canonical $\sigma-$field and a probability measure $\IP$ of the form $ \IP := \mu^{\otimes \Z}$, where $\mu$ is a probability on $[0,1]$. We can now define a probability $\tildeP$ on the space of trajectories, called the \textit{annealed law}, by :
$$ \tildeP = \int_{\Omega} \tildeP_{\omega} d \IP. $$

Introducing the notation $\rho(i) := \frac{1- \omega(i)}{\omega(i)}$, we make the following assumptions on $\mu$ :
\begin{enumerate}
\item[(i)]{$ \mu( \omega(0) = 0) = \mu( \omega(0) = 1) = 0$}
\item[(ii)]{$\E_{\mu}(\log \rho(0)) = 0$}
\item[(iii)]{$\E_{\mu}\left[(\log \rho(0))^2 \right]< \infty$}
\end{enumerate}

A random walk in random environment satisfying the assumptions above is usually called a Sinai walk, referring to the famous article of Sinai \cite{SIN} proving the convergence in law of $X_n / (\log n)^2$ under the annealed law.
\begin{itemize}
\item{Assumption of ellipticity (i) is an irreducibility assumption.}
\item{Assumption (ii) ensures that $\tildeP_{\omega}$ is recurrent $\IP$-almost surely (see \cite{ZEI} for a survey on one dimensional random walks in random environments).}
\item{Assumption (iii) makes it possible to apply Donsker's principle to the potential $V_{\omega}$ (see (\ref{V}) below).}
\end{itemize}

Let us now explain the construction of the internal diffusion limited aggregation cluster for a given environment $\omega \in [0,1]^{\Z}$.

Let $(X^j(n))_{n\in \N}$ be an i.i.d. family of random walks such that $\forall j \in \N$, $(X^j)$ has law $ \tildeP_{\omega}$. We now define our cluster $A(n)$ as a classical internal diffusion limited aggregation cluster using this family of random walks. Define $A(n)$ and the stopping times $(\theta_k)_{k \in \mathbb{N}}$ recursively in the following way :
\begin{eqnarray}
\theta_0 & = & 0 \text{,}\nonumber\\
A(0)     & = & \{ 0 \} = \{ X^0(\theta_0) \} \text{, and for all $j >0$,}\nonumber\\
   \text{ } \theta_j & = & \inf \{ n \geq 0 : X^j(n) \not \in A(j-1) \} \text{,}\nonumber\\
               A(j)  & = & A(j-1) \cup \{ X^j(\theta_j) \} \text{.}\nonumber
\end{eqnarray}
\begin{rem}
We stress the fact that under the annealed law, the same environment is used throughout the construction of $A(n)$.
\end{rem}
As the construction of the cluster uses an i.i.d sequence of walks, the following notations will be helpful to state our theorems,
\begin{equation*}
 P_{\omega}=\tildeP_\omega^{\otimes \N} \textrm{ and }P=\int_{\Omega} P_{\omega} d \IP.
\end{equation*}

The crucial tool introduced by Sinai in \cite{SIN} is the potential associated to a given environment :
 \begin{eqnarray}
V_{\omega}(0) &=& 0 \nonumber \\
V_{\omega}(i) &=& \sum_{k=1}^{i} \ln \rho(k) \text{ if }i \geq 1, \label{V} \\
V_{\omega}(i) &=& -\sum_{k=i+1}^{0} \ln\rho(k) \text{ if }i\leq -1. \nonumber
\end{eqnarray}

For all $n\geq 0$, we define the renormalized potential
\begin{equation*}
 V_{\omega}^{(n)}(t)=\frac{1}{\sqrt{n}}V_{\omega}(\lfloor nt \rfloor).
\end{equation*}

It follows from Donsker's principle and Assumption (iii) that, as $n$ goes to infinity, $(V^{(n)}_{\omega}(t))_{t\in \R}$ converges in law to the Wiener law. Let us notice finally that, by construction,  the cluster $A(n)$ is an interval we will
denote by $[g_n, d_n]$.
We can now state our main result:
\begin{theorem}
\label{th:mainth}
 Under $P$, $d_n/n$ converges in law to the Arcsine law. Namely, for all $0\leq a \leq b \leq 1$,
\begin{equation*}
 P\left(\frac{d_n}{n}\in (a,b)\right)\xrightarrow[n\to\8]{}\int_a^b\frac{1}{\pi\sqrt{x(1-x)}}dx.
\end{equation*}
\end{theorem}

Moreover the following theorem describes the almost sure behavior of the DLA cluster in a typical environment :

\begin{theorem}
 \label{th:quenched}
$\IP-$a.s., with $P_\omega$ probability one, 
\begin{align*}
 \limsup_{n\to\8}\frac{d_n}{n}&=1 \text{, and}\\
 \liminf_{n\to\8}\frac{d_n}{n}&=0.
\end{align*}
\end{theorem}

\section{Proof of Theorem \ref{th:mainth}}

The proof of Theorem \ref{th:mainth} can be decomposed in the three following steps.

\subsection{Good environments}

For each $n$, we define the set of \textit{good environments} which will turn out to be of high probability and on which we will be able to control the position of the cluster at step $n$.

For all cadlag functions $v:\R\to\R$, we define,

\begin{eqnarray}
T_y^+(v) & =&  \inf \left\{t \geq 0,\text{ such that }v(t) \geq y \right\}, \nonumber\\
T_y^-(v) & =& -\sup \left\{t \leq 0,\text{ such that }v(t) \geq y \right\}, \nonumber\\
\bar{y}  &=& \sup \{y \geq 0, \text{ such that } T_y^+ + T_y^- \leq 1 \}.
\end{eqnarray}

The excursion $\alpha$ below the maximum at $T^+_{\bar{y}}$ is defined by
\begin{equation*}
 \alpha = \inf \left\{t \geq 0, \text{ such that }v(T_{\yb}^+ +t) \geq v(T_{\yb}^+)\right\},
\end{equation*}
and the one to the left of $-T_{\yb}^-$ is defined analogously:

\begin{equation*}
 \beta  = \inf \left\{t \geq 0, \text{ such that }v(T_{\yb}^- - t )\geq v(T_{\yb}^-)\right\}.
\end{equation*}

\begin{figure}[ht]
\centering
\includegraphics[scale=0.65]{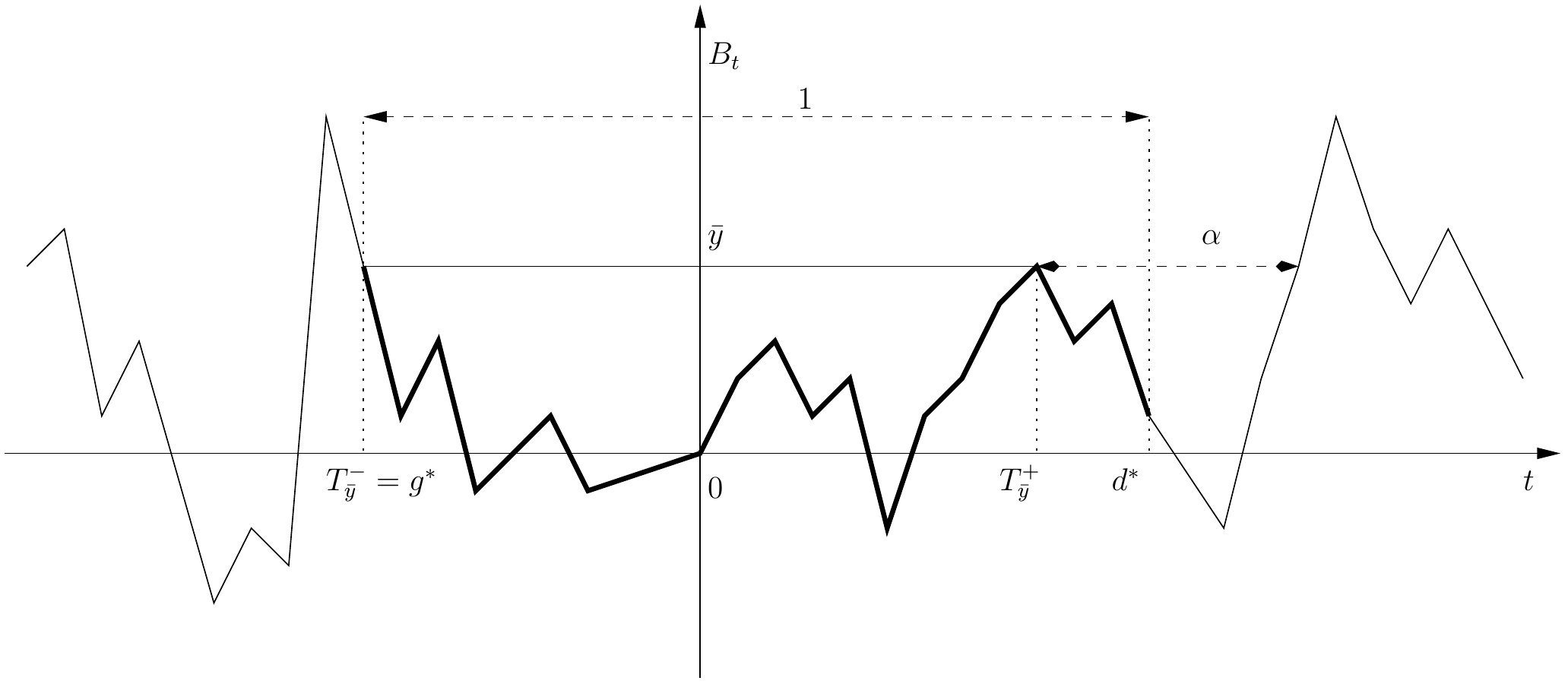}
\caption{\textit{On this example $\alpha>0$ while $\beta=0$. The bold part of the path corresponds to the ``theoretical'' part of the potential that is explored by the cluster.}}\label{fig:dessin1}
\end{figure}

In order to make the computation more readable, we will use for all $n\in\N$ the following notations:

$$ T^{+}_{\yb_n}=T_{\yb}^+(V^{(n)}), \text{    } T^{-}_{\yb_n}=T_{\yb}^-(V^{(n)}), \text{    } \yb_n=\yb(V^{(n)}),$$
$$ \alpha_n=\alpha(V^{(n)}), \text{    } \beta_n=\beta(V^{(n)}).$$
Let $(B_t)_{t\in\R}$ be a real standard Brownian motion defined on an abstract probability space $(\bf{\Omega},\textbf{F},\bP)$. We define with a slight abuse of notation,

$$ T^{+}_{\yb}=T_{\yb}^+(B), \text{    }  T^{-}_{\yb}=T_{\yb}^-(B), \text{    }  \yb=\yb(B), $$
$$ \alpha=\alpha(B), \text{    }  \beta=\beta(B)$$
These notations (as well as $d^*$ and $g^*$ introduced in the next section) are illustrated in Figure \ref{fig:dessin1} in the case of the Brownian motion.
\begin{rem}
Throughout the paper we make use of Donsker's principle for the functionals $\alpha, \beta, T^{+}_{\yb}, T^{-}_{\yb}, \yb$, as well as $d^{*}$ (introduced in the next section), which are not continuous with respect to the Skorohod topology. However, we observe, for all these functionals,  that the set on which they are discontinuous is a subset of trajectories having two local maxima at the same height which has Wiener measure $0$.
\end{rem}

For all $\epsilon>0$ and $n>0$, we define the following events:
\begin{align*}
 B_n^{\epsilon, +} &= \left\{\sup_{\left[-T_{\yb_n}^- - \beta_n -\epsilon ,-T_{\yb_n}^- - \beta_n \right]} V^{(n)}_\omega \geq \yb_n + \frac{n^{1/3}}{\sqrt{n}}\right\} \cap \left\{ \beta_n  < \epsilon \right\},\\
 B_n^{\epsilon, -} &= \left\{\sup_{\left[T_{\yb_n}^+ +\alpha_n,T_{\yb_n}^+ +\alpha_n +\epsilon  \right]} V^{(n)}_\omega \geq \yb_n + \frac{n^{1/3}}{\sqrt{n}}\right\} \cap \left\{ \alpha_n < \epsilon \right\},\\
C_n^{\epsilon, +} &=  \left\{ \sup_{\left[-T_{\yb_n}^-  + \epsilon ,0 \right]} V^{(n)}_\omega \leq \yb_n - \frac{n^{1/3}}{\sqrt{n}}\right\},\\
C_n^{\epsilon, -} &=  \left\{ \sup_{\left[0,T_{\yb_n}^+  - \epsilon  \right]}  V^{(n)}_\omega \leq \yb_n - \frac{n^{1/3}}{\sqrt{n}}\right\}.
\end{align*}
These events can be described as follows. On $B_n^{\epsilon, +}$, the length of the excursion below the supremum to the left of $-T^-_{\yb_n}$ is smaller than $\epsilon$. Furthermore, to the left of this excursion the potential increases enough to build an obstacle for the walk. On $C_n^{\epsilon, +}$, there is no significant obstacle between $0$ and $-T^-_{\yb_n}+\epsilon$. On the intersection of these two events, we expect the left border of the cluster to be close to $-T^-_{\yb_n}$ while the right border should go beyond $T^+_{\yb_n}$. The events indexed by $``-$'' refer to the symmetric situation and can be described in the same way.

\begin{lemma}
\label{lemma:evenement} For all $\epsilon>0$,
 \begin{equation*}
 \lim_{n \rightarrow \infty} \IP\left(\left( B_n^{\epsilon, +} \cap C_n^{\epsilon, +} \right) \cup \left( B_n^{\epsilon, -} \cap C_n^{\epsilon, -} \right) \right) = 1
 \end{equation*}
\end{lemma}
\begin{proof}[Proof of Lemma \ref{lemma:evenement}]
We define the following filtration :
\begin{eqnarray}
\mathcal{F}_0 & = & \sigma \left( V_{\omega}(i) , i \geq 0 \right) \\
\mathcal{F}_k & = & \sigma \left( V_{\omega}(i) , i \geq -k \right) \text{ for all $k>0$.}
\end{eqnarray}

The process $ \left( V_{\omega} ( -i) \right)_{i \geq 0} $ is Markovian and adapted to $\mathcal{F}:=\left( \mathcal{F}_j \right)_{j \geq 0} $. The Markov property at the $\mathcal{F}$-stopping time $(-n T^-_{\yb_n} - n \beta_n)$ and classical properties of random walks yield
$$ \IP \left(\sup_{\left[-T_{\yb_n}^- - \beta_n -\epsilon ,-T_{\yb_n}^- - \beta_n \right]} V^{(n)}_\omega \geq \yb_n + n^{1/3}/\sqrt{n}\right) \xrightarrow[n\to\8]{} 1.$$

The same argument holds for the symmetric case.

It follows from Donsker's principle that
$$ \IP \left( \left\{ \beta_n  < \epsilon \right\} \cup \left\{ \alpha_n  < \epsilon \right\} \right) \rightarrow \bP \left( \left\{ \beta  < \epsilon \right\} \cup \left\{ \alpha  < \epsilon \right\} \right) .$$

As $\left( T_{\yb}^+ \right)_{y \geq 0} $ and $\left( T_{\yb}^- \right)_{y \geq 0} $ are two strictly increasing subordinators without drift, their sum is also a strictly increasing subordinator without drift, and hits (go through) the level $1$ while jumping (see Proposition $1.9$ in \cite{BER}). Furthermore, they are independent so they never jump at the same time (for background on subordinators, we refer to \cite{BER}). As a consequence, under Wiener measure, either $\alpha >0$ and $\beta =0$ or $\alpha =0$ and $\beta >0$, hence
$$ \bP \left( \left\{ \beta  < \epsilon \right\} \cup \left\{ \alpha  < \epsilon \right\} \right) \rightarrow 1. $$

It follows from the same argument that, a.s., one of $T_{\yb}^+$ and $T_{\yb}^-$ is a local maximum (while the other is not). Notice also that $\left(-n^{1/3}/\sqrt{n}, + \infty \right) $ increases to  $\left(0, + \infty \right) $, hence from Donsker's principle,
\begin{equation*}
 \IP \left( C_n^{\epsilon, -} \cap C_n^{\epsilon, +} \right) \xrightarrow[n\to\8]{} \textbf{P} \left( \sup_{[-T_{\yb}^- + \epsilon, T_{\yb}^+ - \epsilon]} B \leq \yb \right).
\end{equation*}
This happens almost surely because the Brownian motion already has a local maximum with value $\yb$ at $T_{\yb}^+$ or $T_{\yb}^-$. 

\end{proof}

\subsection{The quenched localization of $d_n$}

We can now limit our study to the good environments, where with high probability $d_n$ is localized near its theoretical position $d^{*}_n$ that is a deterministic functional of the potential. More precisely, we define for any cadlag function
\begin{equation*}
 d^*=T^+_{\yb}+\ind_{\alpha>\beta}(1-(T^+_{\yb}+T^-_{\yb})),
\end{equation*}
and we will use the following notations :
$$ d^{*}_n=d^{*}(V^{(n)}), \text{    } d^{*}=d^{*}(B),\text{    } g^{*}=1-d^{*}.$$

We will also use the notation
\begin{equation*}
g_n^* = d_n^* - 1.
\end{equation*}

\begin{proposition}
 \label{prop:quenched}For all $\epsilon>0$ and $\eta>0$,
\begin{equation*}
\IP \left( P_{\omega} \left( \left|\frac{d_n}{n}- d_n^*\right| > \epsilon \right) > \eta \right) \xrightarrow[n\to\8]{} 0.
\end{equation*}
\end{proposition}
\begin{proof}[Proof of Proposition \ref{prop:quenched}]

If $(X_n)_{n\geq 0}$ is a Markov chain on $\Z$, we will use, for any $q$ in $\Z$, the notation $\sigma_q$ to denote the hitting time of $q$ by $X_n$, namely
\begin{equation*}
 \sigma_q=\inf\{n\geq 0, X_n=q\}.
\end{equation*}
We recall (see for example \cite{ZEI}) that for any $\omega \in \Omega$,
\begin{equation}
\label{eq:formule}
 \tildeP_{\omega} (\sigma_{-b} < \sigma_a) = \frac{\sum_{i=0}^{a-1} \exp(V_\omega(i) )}{\sum_{i=-b}^{a-1} \exp(V_\omega(i) )}.
\end{equation}
Let $n$ be in $\N\setminus\{0\}$ and $\omega$ be in $B_n^{\epsilon,+}\cap C_n^{\epsilon,+}$, then
\begin{equation}
\label{eq:gauche}
 P_{\omega} (g_n< -n T_{\yb_n}^- -2 \epsilon n)\leq n^2 \exp(- n^{1/3}).
\end{equation}
Indeed suppose that $V_{\omega}(T_{\yb_n}^+)=\yb_n$, then

\begin{align*}
 \tildeP_{\omega} (\sigma_{-n T_{\yb_n}^- -2 \epsilon n} < \sigma_{n \wedge (nT_{\yb_n}^+ + n\alpha_n-1)})  &= \frac{\sum_{i=0}^{n \wedge (nT_{\yb_n}^+ + n\alpha_n-1)-1} \exp(V_\omega(i) )}{\sum_{i=-n T_{\yb_n}^- -2 \epsilon n}^{n \wedge (nT_{\yb_n}^+ + n\alpha_n-1)-1} \exp(V_\omega(i) )}\\
&\leq \frac{n\exp(\yb_n)}{\exp(\yb_n + n^{1/3})} = n \exp(- n^{1/3})
\end{align*}
Now as $n (T_{\yb_n}^+ + T^-_{\yb_n}+2\epsilon+\alpha_n)> n$, the probability in (\ref{eq:gauche}) can be controlled by the probability that one of $n$ independent random walks in the environment $\omega$ exits of the interval $[-n T_{\yb_n}^- -2 \epsilon n , n T_{\yb_n}^+ + n\alpha_n-1 ]$ by the left side, namely
\begin{equation}
\label{eq:control1}
 P_{\omega} (g_n< -nT_{\yb_n}^- -2 \epsilon n)< 1-(1-n\exp(-n^{1/3}))^n \leq n^2 \exp(- n^{1/3})
\end{equation}

Suppose now that $V_{\omega}(T_{\yb_n}^+) > \yb_n$ then
\begin{align*}
  \tildeP_{\omega} (\sigma_{-n T_{\yb_n}^- -2 \epsilon n}<\sigma_{nT_{\yb_n}^+-1}) &=\frac{\sum_{i=0}^{nT_{\yb_n}^+-2} \exp(V_\omega(i) )}{\sum_{i=-n T_{\yb_n}^- -2 \epsilon n}^{nT_{\yb_n}^+-2} \exp(V_\omega(i))}\\
&\leq \frac{n\exp(\yb_n)}{\exp(\yb_n + n^{1/3})} = n \exp(- n^{1/3})
\end{align*}
It follows from the definition of $\yb_n$ that  $n(T_{\yb_n}^+ + T_{\yb_n}^- + \beta_n)\geq n$ which leads to the same kind of control as in (\ref{eq:control1}) and concludes the proof of $(\ref{eq:gauche})$.
We now complete the study of $g_n$ on $B_n^{\epsilon,+}\cap C_n^{\epsilon,+}$ by proving the following inequality,
\begin{equation}
\label{eq:droite}
 P_{\omega} (g_n \geq nT_{\yb_n}^- + \epsilon n)\leq n^2 \exp(- n^{1/3}).
\end{equation}
Using  (\ref{eq:formule}) again, we get
\begin{align*}
 \tildeP_{\omega} (\sigma_{nT_{\yb_n}^+}<\sigma_{-nT_{\yb_n}^- + \epsilon n}) & =  \frac{\sum_{i=0}^{-nT_{\yb_n}^- + \epsilon n-1} \exp(V_\omega(i) )}{\sum_{i=nT_{\yb_n}^+}^{-nT_{\yb_n}^- + \epsilon n-1} \exp(V_\omega(i) )}\\
 &\leq \frac{n\exp (\yb_n - n^{1/3})}{\exp(\yb_n)}=n\exp(-n^{1/3}).
\end{align*}
As $nT_{\yb_n}^+ + nT_{\yb_n}^- + \epsilon n\geq n$, we control the probability of the complementary event in $(\ref{eq:droite})$ with the probability that $n$ independent random walks in the environment $\omega$ exit $[-nT_{\yb_n}^- + \epsilon n, nT_{\yb_n}^+ ]$ through the left side, namely 
\begin{equation*}
 P_{\omega} (g_n \geq -nT_{\yb_n}^- + \epsilon n)\leq 1-(1-n\exp(-n^{1/3}))^n\leq n^2\exp(-n^{1/3}).
\end{equation*}

Notice now that on $B_n^{\epsilon,+}\cap C_n^{\epsilon,+}$, $|T_{\yb_n}^+ - g_n^*|<\epsilon$. With a similar study of the event $B_n^{\epsilon,-}\cap C_n^{\epsilon,-}$ and Lemma \ref{lemma:evenement}, it is easy to complete the proof of Proposition \ref{prop:quenched}.
\end{proof}

\subsection{Characterization of the law of $d^*$}
Using Donsker's principle and Proposition \ref{prop:quenched}, we conclude that $d_n/n$ converges in law towards $d^{*}$. To complete the proof of Theorem \ref{th:mainth}, we characterize the law of $d^{*}$.

\begin{lemma}
\label{lemma:detoile}
 The law of $d^*$ is the Arcsine law.
\end{lemma}
\begin{proof}[Proof of Lemma \ref{lemma:detoile}]
We note $(L_t)_{t\geq 0}$ the local time at $0$ of $B$ and
\begin{equation*}
 \tau_t=\sup\{u\geq 0,\ L_u<t\}
\end{equation*}
the left continuous inverse of $L$.
\begin{remark}
\label{rem:gauche}
The process $(T_y^+)_{y\leq0}$ is left continuous. Consequently, even if it is
unusual, we prefer to work with the left continuous version of the inverse local time of
$B_t$. We will also use the name subordinator for a left continuous increasing process
with independent increments.
\end{remark}
For all $t\geq 0$, denote 
\begin{equation*}
 A^+_t=\int_0^t 1_{B_s>0} ds,
\end{equation*}
the time spent by $B$ above $0$ up to time $t$ ; and similarly
\begin{equation*}
 A^-_t=\int_0^t 1_{B_s<0} ds.
\end{equation*}
Now $A^+_{\tau_t} + A^-_{\tau_t}=\tau_t$ is a subordinator and, as $\tau_t$ follows the same law as the hitting time of $t$ by the Brownian motion, it is a $1/2-$stable subordinator. Moreover it is a consequence of Ito's excursion theorem (see Theorem 2.4 of \cite{RY}) that both processes $(A^+_{\tau_t})$ and $(A^-_{\tau_t})$ are independent subordinators with the same law. Each of them deals indeed respectively with one of the disjoint set of positive and negative excursions of $B$ up to time $\tau_t$. Hence 
\begin{equation}
\label{eq:keyidentity}
 (A^+_{\tau_t},A^-_{\tau_t})\stackrel{(law)}{=}\frac{1}{4}(T^+_t,T^-_t)
\end{equation}
(we refer to the lecture notes of Marc Yor \cite{Y_CAR} for explanations and applications of this
identity).
Hence
\begin{align*}
 \bar{y}&\eql \sup\{t\geq 0,\ A^+_{\tau_t}+A^-_{\tau_t}\leq 1\}\\
	&\eql \sup\{t\geq 0,\ \tau_t\leq 1\}\\
	&\eql L_1.
\end{align*}
From the choice of the left continuous version of the inverse local time (see Remark \ref{rem:gauche}) we deduce
\begin{equation}
\label{eq:loiyb}
 T^+_{\bar{y}}\eql A^+_{\tau_{L_1}}=A^+_{g_1}
\end{equation}
where
\begin{equation*}
 g_1=\sup\{t<1,\ B_t=0\}
\end{equation*}
is the last zero of $(B_t)_{t\geq 0}$ before $1$.
Let us remind that
\begin{equation}
\label{eq:defdetoile}
 d^*=T_{\bar{y}}^+ + 1_{\alpha>\beta} (1-(T^+_{\bar{y}}+T^-_{\bar{y}})).
\end{equation}
Notice now that $(T^+_{\bar{y}},T^-_{\bar{y}})$ is independent of $1_{\alpha>\beta}$. Indeed, using $(\ref{eq:keyidentity})$, it is equivalent to check the independence of the sign of $B_1$ and $(A^+_{g_1},A_{g_1}^-)$. Gathering this independence with $(\ref{eq:defdetoile})$ and $(\ref{eq:loiyb})$, we obtain
\begin{equation}
\label{eq:horiz}
 d^*\eql A^+_{g_1}+\epsilon (1-(A^+_{g_1}+A^-_{g_1}))=A^+_{g_1}+\epsilon (1-g_1),
\end{equation}
where $\epsilon$ is a Bernoulli variable with parameter $1/2$ independent of all other variables. The decomposition of the path of the Brownian motion on $[0,1]$ into a path on $[0,g_1]$ and an incomplete excursion of independent sign on $[g_1,1]$ yields
\begin{equation}
\label{eq:verti}
 A^+_{g_1}+\epsilon (1-g_1)\eql A^+_1.
\end{equation}
We conclude by recalling Paul Levy's well known result (\cite{LEVY})which states that $A^+_1$ follows the Arcsine law.
\end{proof}

\begin{remark}
 Notice that the Proof of Lemma $\ref{lemma:detoile}$ has a non-trajectorial nature as shown by the key identity in law $(\ref{eq:keyidentity})$, which allows to reduce the study of the functional of two independent Brownian motion to that of the functional of a single one.
\end{remark}
\begin{remark}
 The random sign $\epsilon$ in identity $(\ref{eq:horiz})$ corresponds to the left-right symmetry of the problem, whereas in the classical decomposition $(\ref{eq:verti})$ it has an up-down meaning.
\end{remark}


\section{Proof of Theorem \ref{th:quenched}}

We will only prove the first statement as the second one can be easily deduced using the symmetry of the model.
\begin{lemma}
 \label{lemma:droite}
$\IP$-p.s.,
\begin{equation*}
 \{T_{\yb_n}^+>1-\epsilon\}\qquad\textrm{i.o.}
\end{equation*}
\end{lemma}
\begin{proof}[Proof of Lemma \ref{lemma:droite}]
Let $k$ be in $\mathbb{N}^*$. It is a well known result that there exist $k$ independent Brownian motions $(B^1,\cdots,B^k)$ on $(\bf{\Omega},\textbf{F},\bP)$ such that
\begin{equation*}
 (V^{(n)},\cdots,V^{(n^k)})\stackrel{(law)}{\Rightarrow}(B^1,\cdots,B^k).
\end{equation*}
We will also use the fact that for any real Brownian motion $B$, $$\bP(T^+_{\yb}(B)\geq 1-\epsilon)>0,$$ see for example $(\ref{eq:loiyb})$. Now
\begin{equation*}
 \liminf\{T_{\yb_n}^+\leq 1-\epsilon \} \subset \liminf\left\{ \{T_{\yb_n}^+\leq 1-\epsilon \}\cap\cdots\cap \{T_{\yb_n^k}^+\leq 1-\epsilon \} \right\}
\end{equation*}
and
\begin{align*}
 \IP(\liminf\{T_{\yb_n}^+\leq 1-\epsilon \}) &\leq \liminf \bP \left( \left\{ T_{\yb_n}^+   \leq 1-\epsilon \right\}\cap    \cdots \cap              \left\{T_{\yb_n^k}^+ \leq 1-\epsilon    \right\} \right) \\
					                                   &\leq         \bP \left(         T_{\yb}^+(B^1)\leq 1-\epsilon         \right) \cdots \bP                \left( T_{\yb}^+(B^k)\leq 1-\epsilon           \right) \\
                                             &\leq         \bP \left(         T_{\yb}^+(B^1)\leq 1-\epsilon         \right)^k	
\end{align*}
As $k$ can be chosen arbitrarily big, this concludes the proof of Lemma \ref{lemma:droite}.
\end{proof}
We deduce from Lemma \ref{lemma:droite} and the previous study of $C_n^{\epsilon,-}$ (see the Proof of Theorem $\ref{th:mainth}$) that $C_n^{\epsilon,-}\cap \{ T_{\yb_n}^+>1-\epsilon\}$ occurs infinitely often. Fix $\omega$ in $\Omega$ and $n$ in $\mathbb{N}$ such that $\omega\in C_n^{\epsilon,-}\cap \{ T_{\yb_n}^+>1-\epsilon\}$. Formula $(\ref{eq:formule})$ yields
\begin{align*}
 \tildeP_{\omega}(\sigma_{nT_{\yb_n}^-}<\sigma_{n(T_{\yb_n}^+-\epsilon)})& \leq \frac{\sum_{i=0}^{n(T_{\yb_n}-\epsilon)-1}\exp(V_\omega(i))}{\sum_{i=nT_{\yb_n}^-}^{n(T_{\yb_n}-\epsilon)-1}\exp(V_\omega(i))}\\
 &\leq \frac{n e^{\yb_n-n^{1/3}}} {e^{\yb_n}} \leq n e^{-n^{1/3}}
\end{align*}
The probability of the event $\{d_n\leq T_{\yb_n}^+-\epsilon\}$ is smaller than the probability that one of $n$ independent random walks in the environment $\omega$ exits $[-nT_{\yb_n}^+,n(T_{\yb_n}^+-\epsilon)]$ through the left side. Hence,
\begin{align*}
 P_{\omega}(d_n\leq T^+_{\yb_n}-\epsilon)& \leq 1-(1-ne^{-n^{1/3}})^n\\
					& \leq n^2e^{-n^{1/3}}.
\end{align*}
We conclude using Borel Cantelli's Lemma on a subsequence $(n_j)_{j\geq 0}$ such that $C_ {n_j}^{\epsilon,+}\cap \{T_{\yb_n}^+>1-\epsilon\}$ holds for all $j\geq 0$ (the $\IP$-a.s. existence of such a subsequence is a consequence of Lemma \ref{lemma:droite}).
\begin{flushright}
$\square$
\end{flushright}

\textbf{Acknowledgements}

It is a pleasure to thank Professor Marc Yor for pointing out to us the key identity (\ref{eq:keyidentity}).

\bibliographystyle{plain}
\bibliography{sinai}
\end{document}